\documentclass{amsart}
\usepackage{amsmath,amsthm,amssymb,xypic}
\usepackage[all]{xy}

\title{Equivalent birational embeddings
  II: divisors}

\author{Massimiliano Mella
\and Elena Polastri}
\address{Dipartimento di Matematica\\
Universit\`a di Ferrara\\
Via Machiavelli 35\\
44100 Ferrara, Italia} \email{mll@unife.it \ plslne@unife.it}

\date{November 2010}
\thanks{Partially supported by Progetto PRIN 2008 ``Geometria
  sulle variet\`a algebriche'' MUR}

\theoremstyle{plain}

\newtheorem{theorem}{Theorem}[section]

\newtheorem{proposition}[theorem]{Proposition}
\newtheorem{lemma}[theorem]{Lemma}
\newtheorem{corollary}[theorem]{Corollary}

\newtheorem{example}[theorem]{Example}

\theoremstyle{definition}

\newtheorem{definition}[theorem]{Definition}

\theoremstyle{remark}

\newtheorem*{claim}{Claim}

\newtheorem{case}[theorem]{Case}
\newtheorem{con}[theorem]{Construction}

\newtheorem{remark}[theorem]{Remark}

\DeclareMathOperator{\rk}{rk} 
\DeclareMathOperator{\Bsl}{Bs}

\DeclareMathOperator{\mult}{mult}

\DeclareMathOperator{\Sing}{Sing}

\DeclareMathOperator{\elm}{elm}

\newcommand{\QED}{\ifhmode\unskip\nobreak\fi\quad {\rm Q.E.D.}} 

\newcommand\Ka{\overline{\kappa}}
\newcommand\map{\dasharrow}
\newcommand\iso{\cong}

\newcommand{\f}{\varphi}
\newcommand{\A}{\mathcal{A}}
\newcommand{\B}{\mathcal{B}}
\newcommand{\C}{\mathbb{C}}

\newcommand{\E}{\mathcal{E}}
\newcommand{\F}{\mathbb{F}}
\newcommand{\G}{\mathcal{G}}

\newcommand{\I}{\mathcal{I}}

\renewcommand{\L}{\mathcal{L}}

\newcommand{\N}{\mathbb{N}}
\renewcommand{\O}{\mathcal{O}}
\renewcommand{\P}{\mathbb{P}}
\newcommand{\Proj}{\mathbb{P}}
\newcommand{\Q}{\mathbb{Q}}

\newcommand{\s}{\mathbb{S}}

\newcommand{\rat}{\dasharrow}

\begin{document}

\begin{abstract} Two divisors in $\P^n$
are said to be Cremona equivalent if
there is a Cremona modification sending
one to the other. We produce infinitely
many non equivalent divisorial
embeddings of any variety of dimension
at most $14$. Then we study the special
case of plane curves and rational
surfaces. For the latter we
characterise surfaces Cremona
equivalent to a plane. \end{abstract}

\maketitle

\section*{Introduction}
Let $X_1, X_2\subset\P^n$ be two birationally equivalent
projective varieties.  It is natural to ask if there exists  a
Cremona transformation of $\P^n$ that maps $X_1$ to $X_2$, in this
case we say that $X_1$ and $X_2$ are Cremona equivalent, see
Definition \ref{def:cremona} for the precise statement. This is
somewhat related to the  Abhyankar--Moh problem, \cite{AM} and
\cite{Je}. Quite surprisingly the main theorem in \cite{MP} states
that this is the case as long as the codimension of $X_i$ is at
least 2. In this work we want to study the case of divisors. It is
easy to give examples of pairs of birationally equivalent divisors
that are not Cremona equivalent. Our approach is to use Log
Minimal Model Program, LMMP for short, and its variants like
Sarkisov Theory and Noether Fano inequalities. Via these
techniques we are able to produce many examples of these
inequivalent embeddings in arbitrary dimensions. Moreover we prove
that any irreducible and reduced variety, of dimension at most
$14$, admits infinitely many non Cremona equivalent divisorial
embeddings. This shows the difficulty of classifying inequivalent
embeddings. This is why we next concentrate on two special classes
of divisors: plane curves and rational hypersurfaces. One can look
at Cremona equivalence as the action of the Cremona group of
$\P^n$ on the Hilbert scheme of divisors.  For $\P^2$ both Cremona
group and  Sarkisov theory are well understood. This allows us to
address a question posed in \cite{Ii} about
minimal degree curves. In Theorem
\ref{th:main} we give a necessary and
sufficient condition for a curve to be of
minimal degree under the Cremona
equivalence. Even if this
description is not fully satisfactory,
it is probably the best one. Pairs 
like those described in Example \ref{ex:bad} show that
nested infinitely near 
singularities give unpredictable
behaviour  with respect to the
Cremona action, and
a partial resolution of the
singularities is needed to decide
minimality. The second class of divisors
we consider is that 
of rational hypersurfaces. A nice result proposed by Coolidge
state that a plane curve is Cremona equivalent to a line if and
only if $\Ka(\P^2,C)<0$, see
  Definition \ref{def:lkod}. This statement has been proved and
somewhat strengthened by Kumar and Murthy, \cite{KM}. Our approach
with LMMP techniques gives a new proof of it and suggests a
possible extension  in arbitrary dimension. The idea is to
consider a log resolution, say $(S,\overline{C})$,  of the pair
$(\P^2, C)$ and translate the hypothesis on Kodaira dimension into
a geometric restriction to the possible contractions occurring
along a LMMP directed by $\overline{C}$. In this way we end up on
log varieties we are able to control.
 This reminded us the
$\sharp$-MMP, \cite{Me},  where again
 numerical constrains where used to
 control the birational modification
 occurring along a LMMP. With these two
 constructions in mind we are able to
prove a Coolidge type statement also for
rational surfaces in $\P^3$, Theorem \ref{th:3c}.

\smallskip

{\it Acknowledgement}. We are much indebted to Alberto Calabri and
Ciro Ciliberto for discussions and suggestions. Especially from
preventing us to cruise in a wrong direction in the definition of
minimal degree curves. Part of this work
has been done during a pleasant stay at
the MSRI. Massimiliano 
Mella would like to thank the Institute 
for providing a perfect environment. 
The anonymous referee has also to be warmly
thanked for a careful reading and many
suggestions that improved exposition and
readability.

After this paper was completed Ciro
Ciliberto brought to our attention his work
with Alberto Calabri, \cite{CC}, devoted
to a detailed study of minimal degree plane
curves. 

\section{Notations and preliminaries}
We work over an algebraically closed field of characteristic zero.
We are inte\-re\-sted in birational transformations of log pairs.
For this we introduce the following definition.

\begin{definition}\label{def:cremona} Let $D\subset X$ be
  an irreducible and  reduced divisor on
  a normal variety $X$.
We say that $(X,\alpha D)$ is birational
  to $(X',\alpha D')$, for  $\alpha\in\Q$, if there exists a
  birational map $\f:X\map X'$ with
  $\f_*(D)=D'$.
Let $D, D'\subset\P^n$ be
   irreducible reduced divisors then we say that $D$ is Cremona
  equivalent to $D'$ if
  $(\P^n,D)$ is birational to $(\P^n,D')$.
\end{definition}

Let us proceed recalling a well known class of singularities.

\begin{definition}
Let $X$ be a normal variety and $D=\sum d_i D_i$ a $\Q$-Weil divisor, with $d_i\leq 1$.
Assume that $(K_X+D)$ is $\Q$-Cartier. Let $f:Y\to X$ be a log
resolution of the pair $(X,D)$ with
$$K_Y=f^*(K_X+D)+\sum a(E_i,X,D) E_i$$
We call
\begin{multline*}
disc(X,D):=\min_{E_i}\{a(E_i,X,D)|\\
\mbox{$E_i$ is an $f$-exceptional
  divisor for some log resolution} \}
\end{multline*}

Then we say that $(X,D)$ is
 \[\begin{array}{ccc}
 \left.\begin{array}{l}
terminal\\
canonical\\
\end{array}
 \right\}&\mbox{if $disc(X,D)$}&
\left\{\begin{array}{l}
>0\\
\geq 0\\
\end{array}
 \right.
\end{array}
\]
\end{definition}

\begin{remark} In the case of smooth surfaces
  the notion of canonical singularities
  has a nice and natural
  interpretation. Any log resolution of
  a smooth surface can be obtained via
  blow up of smooth points. Hence a pair $(S,D)$, with
  $S$ a smooth surface has canonical
  singularities if and only if $\mult_p
  D\leq 1$ for any point $p\in S$. 

Note further that one direction is true
in any dimension. Assume that $X$ is
smooth and  $\mult_p
D\leq 1$ for any point $p\in X$. Let  $f:Y\to X$ be a
smooth blow up, with exceptional divisor
$E$. Then $K_Y=f^*(K_X)+aE$ for some
positive integer $a$ and $a(E,X,D)\geq
a-1\geq 0$. This proves that 
 $(X,D)$ has canonical
singularities if $X$ is smooth and $\mult_pD\leq
1$ for any $p\in X$.
\label{rem:mult}
\end{remark}
To investigate pairs of divisor and
variety it is useful to introduce the
following notation. 
\begin{definition}\label{def:lkod}
 Let $(X,D)$ be a pair, with $D=\sum d_i
 D_i$ and $0\leq d_i\leq 1$. Fix $f:Y\to X$ a log resolution with  $D_Y=f^{-1}_*D$ the strict transform.
We indicate
with $\Ka(X,D):=\kappa(Y,D_Y)$, where
$\kappa(Y,D_Y)$ is the transcendence
degree of of the ring $\oplus_m
H^0(Y,m(K_Y+D_Y))$ minus 1.
\end{definition}

For future reference we recall  a
technical result on pseudoeffective
divisors, i.e. the closure of effective divisors.
\begin{lemma}\label{lem:psef}
  Let $(X,D_X)$ and $(Y,D_Y)$ be
  birational pairs with canonical
  singularities. 
Then $K_X+D_X$ is pseudoeffective if and
  only if $K_Y+D_Y$ is pseudoeffective.
\end{lemma}
\begin{proof}
  It is enough to prove the statement
  when there is a morphism $f:Y\to X$
  such that $D_Y=f^{-1}_*D_X$. Assume
  that $K_Y+D_Y$ is pseudoeffective and
  let $A$ be a general ample divisor on
  $Y$. Then for any $\epsilon>0$
  $K_Y+D_Y+\epsilon A$ is
  effective. Hence $f_*(K_Y+D_Y+\epsilon
  A)=K_X+D_X+f_*(\epsilon A)$ is
  effective, and $K_X+D_X$ is
  pseudoeffective. Assume that $K_X+D_X$
  is pseudoeffective then
  $K_Y+D_Y=f^*(K_X+D_X)+\Delta$, for
  some effective divisor $\Delta$. Hence
  $K_Y+D_Y$ is pseudoeffective.
\end{proof}

While working on plane curves we
frequently use ruled surfaces to fix the
notations we recall the following
definition.

\begin{definition} Let
  $\F_a:=\P(\O\oplus\O(-a))$ be the
  Hirzebruch surface. Then  the exceptional section
  ( a fibre for $a=0$) will be called
  $C_0^a$ and a fibre $f$.
\end{definition}

In the final section we shall use
some standard Cremona maps of $\P^3$, see
for instance \cite{SR}, and birational modifications of
scrolls. We find it
convenient to group them here. Following
\cite{SR} a Cremona map given by forms of
degree $a$ and such that inverse is given by
forms of degree $b$ is said to be of
type $T_{a,b}$.

\begin{con}[$T_{2,3}$]\label{con:23}Let $l\subset\P^3$ be a line and
${p_1,p_2,p_3}$ three general
points. The linear system of quadrics
through this configuration is homaloidal
and gives rise to a Cremona
transformation of type $T_{2,3}$. For
our purpose the following  facts
are important:
\begin{flushleft}
\ the rational normal curves
  secant to $l$ and passing through the
  $p_i$'s are sent to lines;\\
\ the plane $P$ spanned by the
  $p_i$'s is contracted to a line by the
  linear system of conics through the
  $p_i$'s and the point of intersection
  $l\cap P$.
\end{flushleft}
\end{con}
\begin{con}[$T_{3,3}$] Let $\Gamma\subset\P^3$ be a rational
normal curve. Two general cubics
containing $\Gamma$ intersect along
$\Gamma\cup C$. Then the linear system
of cubics through $C$ is homaloidal and
produces a Cremona transformation of
type $T_{3,3}$. It is a pleasant
exercise to check that any rational
cubic with isolated singularities is
contained in one such linear system, and
is therefore Cremona equivalent to a plane.
\label{con:33}
\end{con}

\begin{definition} Let $X:=\P(\E)\to
  \P^1$ be a scroll with $F$ a
  fiber. Fix a point $x\in F\subset
  X$. Then the elementary
  transformation centred at $x$ is
$$\elm_x:X\rat X'$$
the composition of the blow up of $x$
and the contraction of the strict
transform of $F$. Then $X'$ is still a
scroll over $\P^1$.

Let $X:=\P(\E)\to W$ be a scroll over a
surface, and $\Gamma\subset X$ a smooth
curve section. Let $\pi:X\to W$ be the
scroll structure and
$S_F:=\pi^*(\pi_*(\Gamma))$. Then the
elementary transformation centred at
$\Gamma$ is
$$\elm_\Gamma:X\rat X'$$
the composition of the blow up of
$\Gamma$ and the contraction of the
strict transform of $S_F$. Then $X'$ is
still a scroll over $W$.
\end{definition}

\begin{con}\label{con:scroll} A nice feature of these
  elementary transformations is the
  following. Let $X$ be a 3-fold scroll
  over $W$, and $S\subset X$ a smooth
  surface section. Then $\pi_{|S}:S\to
  W$ is a birational map. Let $D\subset
  X$ be a general surface section and
  $\Gamma:=S\cap D$. Then we have
   $\elm_{\Gamma}(S)\iso
  W$.\label{rem:scroll}
\end{con}

\section{Existence results}
We expect that any projective variety has infinitely many non
Cremona equivalent divisorial embeddings. A slight variation of
\cite[Lemma 3.1]{MP} together with the remarkable work of Mather,
\cite{Ma1}, \cite{Ma2}, \cite{Ma3} allows us to prove it under a dimensional
bound.
\begin{theorem} Let $X$ be an
  irreducible reduced projective variety
  of dimension $k$. Assume that $k\leq
  14$ then there are infinitely many non
  Cremona equivalent embeddings of $X$
  in $\P^{k+1}$.
\label{th:Mather}
\end{theorem}

To prove the theorem we recall in a
slightly generalised form \cite[Lemma
3.1]{MP}.

\begin{lemma} Let $X^{n-1}$ be an
  irreducible and reduced projective
  variety. Let $\L$ and $\G$ be linear
  systems of projective dimension
  $n$. Let $\f_\L$ and $\f_\G$ be the
  induced map into $\P^n$. Assume that
  $\f_\L(X)$ and $\f_\G(X)$ are divisors
  of degree, respectively, $l$, and $g$,
  with $l>g$. If $\f_\L(X)$ is
 Cremona equivalent to $\f_\G(X)$ then
 the pair
  $(\P^n,\frac{n+1}l \f_\L(X))$ has not
  canonical singularities.
\label{lem:NF}
\end{lemma}
\begin{proof}
Let $\Phi:\P^n\rat\P^n$ be a Cremona
equivalence   between $\f_\L(X)$ and
$\f_\G(X)$.
Fix a resolution of $\Phi$
\[
\xymatrix{
&Z\ar[dl]_{p}
 \ar[dr]^q &\\
 \P^n\ar@{.>}[rr]^\Phi&&\P^n\\
}
\]

 Then we have
$$\O_Z\sim
p^*(\O(K_{\P^n}+\frac{n+1}l\f_\L(X))=K_Z+\frac{n+1}l X_Z-\sum
a_iE_i$$
and
$$
q^*(\O(K_{\P^n}+\frac{n+1}l\f_\G(X))=K_Z+\frac{n+1}l X_Z-\sum
b_iF_i $$ where $E_i$, respectively $F_i$, are $p$, respectively
$q$, exceptional divisors. Let $r\subset\P^n$ be a general line in
the right hand side $\P^n$. Then the hypothesis $l>g$ yields
$$0 > q^{-1}r\cdot (q^*(\O(K_{\P^n}+\frac{n+1}l\f_\G(X))+\sum
b_iF_i)=(\sum a_iE_i)\cdot q^{-1}r$$
This proves that at least one $a_i<0$ proving
the claim.
\end{proof}

\begin{proof}[Proof of Theorem \ref{th:Mather}]
Without loss of generalities we may
assume that $X$ is smooth of dimension
$k$. 
Fix $\A$, $\B$ very
ample linear system of degree
at least $(k+1)(k+2)$. Consider two general
sublinear system $\L\subset\A$ and
$\G\subset\B$ both of projective
dimension $k+1$. Let
$X_1=\f_\L(X)\subset\P^{k+1}$ and
$X_2=\f_\G(X)\subset\P^{k+1}$. We
already noticed in Remark \ref{rem:mult} 
that if the multiplicities of a pair are
bounded by 1 the pair has canonical
singularities.
Hence 
Lemma \ref{lem:NF} shows that $X_1$ and
$X_2$ are not Cremona equivalent as long
as $\deg X_1>\deg X_2$ and
$$\max_{x\in
  X_1}\{\mult_xX_1\}\leq \frac{\deg
  X_1}{k+2}$$

 If $k=\dim X\leq 14$ then by Mather
 transversality \cite{Ma1} \cite{Ma3}, see
also \cite{EB}, a general projection has
points of multiplicity at most
$k+1$. Then $X_1$ is not Cremona
equivalent to $X_2$.
\end{proof}

\begin{remark} Mather's results are
  optimal as explained by Lazarsfeld,
  \cite[Theorem 7.2.19]{Laz} see also \cite{EB}. In
  general it is not known a bound,
  depending only on the dimension,
  for the singularities of a general
  projection. In \cite{EB} it is proposed a
  conjectural bound  on the Mumford
  regularity instead of
 the multiplicity of the fiber of a
 general projection. 
\end{remark}

Via Lemma \ref{lem:NF} it is easy to
give examples of inequivalent embeddings
in arbitrary dimension.
Consider a smooth codimension two subvariety
$X_{a,b}\subset\P^{k+2}$ given by the complete
intersection of a form of degree $a$
with a form of degree $b$, for $a\leq b$. Then a
general projection to $\P^{k+1}$ has
degree $ab$ and
points of multiplicity at most $a$. On
the other hand projecting from a general
point of $X$ produces a birational embedding
of $X$ of degree $ab-1$. By Lemma
\ref{lem:NF} these two embeddings are
not birational equivalent as long as
$b\geq k+2$. Theorem \ref{th:Mather} and
these examples show that the classification of
inequivalent embeddings  is almost
hopeless for a general variety. On the
other hand there are special classes of
varieties for which something more can
be said.
\section{Plane curves}

In this section we study the Cremona
equivalence for plane curves. Our aim is
to describe minimal degree
representative in each class of Cremona equivalence.


\begin{definition} Let $C\subset S$
  be an irreducible reduced curve on a
  smooth surface $S$.  Let $\Sing^{\infty}(C)=\{p_i\}$
  be the set of, possibly infinitely
  near, points with $m_i=\mult_{p_i}C>1$, and  $\mult(C)=\{m_i\}$ the associated set of
  multiplicities. We always order
  the $p_i$ in such a way that $m_i\geq
  m_{i+1}$.

\end{definition}


It is clear that any pair $(\P^2,C)$ is
birational to a pair $({\mathbb
  F}_a,\tilde C)$. Our first aim is to
choose a nice representative of the pair
$(\P^2,C)$. This is a slight variation
on the usual Sarkisov program for log
surfaces, \cite{BM}.

\begin{proposition}\label{pro:model} Let $C\subset\P^2$ be an
  irreducible reduced curve of degree $d$. Then
  $(\P^2,C)$ is birational
  to one of the following:

\begin{tabular}{|ll|lp{5cm}|}
\hline&&&\\ $(\P^2,l)$& $l\subset \P^2$ a
  line&$(\P^2, \overline{C})$&
 $\overline{C}\sim\O(d')$\\
&&& $K_{\P^2}+\frac3{d'} \overline{C}\equiv 0$
\\&&&$(\P^2,\frac3{d'} \overline{C})$
  terminal\\
\hline&&&\\
 $(\F_0,\overline{C})$&
$\overline{C}\sim \alpha C^0_0+\beta
  f$& $(\F_a,\overline{C})$&$\overline{C}\sim \alpha C^a_0+\beta
  f$\\
& $K_{\F_a}+\frac2\alpha \overline{C}$ nef&&
$K_{\F_a}+\frac2\alpha \overline{C}$ nef\\
&$(\F_0,\frac2\alpha \overline{C})$
  terminal&& $(\F_a,\frac2\alpha \overline{C})$ canonical\\
&&&$(\F_a,\frac2\alpha \overline{C})$
  terminal along $C^a_0$\\&&&\\
\hline
\end{tabular}

\end{proposition}
\begin{proof}We prove the statement by
  induction on the degree of $C$.
Assume that
  $(\P^2,\frac3{d} C)$ is not
  terminal hence
\begin{equation}m_1\geq \frac{d}3
\label{eq:nc}
\end{equation}
If $d-m_1=1$, then
$C$ is Cremona equivalent to a
line.

Assume that $d-m_1\geq 2$ and let
$\epsilon:\F_1\to\P^2$  be the blow
up of $p_1$. Let
$$\overline{C}=\epsilon^{-1}_*C\sim (d-m_1)C^1_0+d f$$
be the strict transform. Then
$$\overline{C}\cdot C^1_0=m_1,$$
and we have
$$(K_{\F_1}+\frac2{d-m_1}\overline{C})\cdot
C^1_0=-1+\frac{2m_1}{d-m_1}\geq 0$$
In
particular $$K_{\F_1}+\frac2{d-m_1}\overline{C}
\ \ \mbox{is nef}$$

In the following we will always denote
$\overline{C}$ the strict transform of
$C$ on the various birational models.
\begin{claim}
There is a chain of elementary
transformations $\f:\F_1\rat\F_b$, such
that $(\F_b,\frac2{d-m_1}\overline{C})$
has 
canonical singularities and it is terminal in a neighbourhood of
$C^b_0$.
\end{claim}
\begin{proof}[Proof of the Claim]
Let $p\in\overline{C}$ be a non
canonical point, then
$\mult_p\overline{C}>\frac{d-m_1}2$. Let
$\phi:\F_1\rat\F_{1\pm1}$ be the
elementary transformation centered in
$p$. Then $\phi(\overline{C})$ acquires
a point, say $q$, of multiplicity $m_q$
with
$$m_q=d-m_1-\mult_p\overline{C}<\frac{d-m_1}2. $$
Hence after finitely many of those
modifications this produces a
 pair
$(\F_h,\overline{C})$ with canonical
 singularities. 
Next we want
to move the canonical singularities away
from the exceptional section $C^h_0$. To
do this we apply again a chain of
elementary transformations centered on canonical points
contained in the exceptional
section. 
Let us understand better the latter
chain. 
Let $p\in \overline{C}\cap C_0$ be a non
terminal point. That is   $\mult_p\overline{C}=\frac{d-m_1}2$. Let $\phi:\F_h\map
\F_{h+1}$ be the elementary
transformation centered at $p$.
Then $\phi(\overline{C})$ has a new
point $q\not\subset C_0^h$ of
multiplicity
$$m_q=d-m_1-\frac{d-m_1}2=\frac{d-m_1}2.$$
At the same time the singularity along
$C_0^{h+1}$ has been simplified.
This means that after finitely many
elementary transformations we obtain a pair
$(\F_b,\frac2{d-m_1}\overline{C})$
canonical and terminal in a neighbourhood of
$C^b_0$.
\end{proof}

Let
$$\overline{C}\sim (d-m_1)C^b_0+\beta f $$
If $b\geq 2$ we have finished because $K_{\F_b}\cdot C^b_0\geq 0$. Assume that $b=1$ and
$$K_{\F_1}+\frac2{d-m_1}\overline{C}
\ \ \mbox{is not nef}$$
then
$$ \frac2{d-m_1}(m_1-d+\beta)<1$$
that is
$$ 2\beta<3(d-m_1)$$
By equation (\ref{eq:nc}) this yields
$$ \beta<d.$$
In other words $C$ is birational to a plane
curve of degree $\beta<d$
and we conclude by induction on the degree.
Assume finally that $b=0$ and
$$K_{\F_0}+\frac2{d-m_1}\overline{C}
\ \ \mbox{is not nef}$$
then $\beta<(d-m_1)$. This time we
change the ruling of $\F_0$ and repeat the
argument for $(\F_0,\frac2\beta
\overline{C})$. The coefficients used are in
$\frac2\N$, therefore after finitely
many steps we find the required model.
\end{proof}

\begin{definition} Let $C\subset\P^2$ be
  a curve. A standard model of
  $(\P^2,C)$ is a birational pair
  $(S,\overline{C})$ obtained via the
  construction of Proposition  \ref{pro:model}.
\end{definition}

The main difficulty in using these
models is
that they are not always
unique.

\begin{remark}\label{rem:final} Thanks to the
  uniqueness of minimal models for
  surfaces the nef models with terminal
  singularities are unique. Note that in presence of
  canonical sin\-gu\-la\-ri\-ties the uniqueness is
  lost.
 Let $C\subset\P^2$ be a 6-ic
  curve with an ordinary double point
  and a tac-node. If we blow up the
  ordinary double point we end up with a
  model $(\F_1,\frac12C)$. While
  resolving the tac-node we produce a
  model
$(\F_2,\frac12C)$.
\end{remark}

Even if not unique the above models
allow us to choose
minimal degree representative in every
Cremona class. 
Let us start with the following.

\begin{proposition}
\label{pro:unique} Let $(S, C)$ and $(S', C')$ be two birational,
not biregular, mo\-dels in the list of \emph{Proposition
\ref{pro:model}}. Then both $S$ and $S'$ are ruled surfaces. Let
\linebreak$\Phi:(S,C)\map (S',C')$ be a birational map of the
pairs and assume that \linebreak$C\sim\alpha C_0+\beta f$ and
$C'\sim \alpha'C^{\prime}_0+\beta' f$. Then
\begin{itemize}
\item[i)] $\alpha=\alpha'$;
\item[ii)] if $\kappa(S,2/\alpha C)=1$, then
  $\Phi$ is an isomorphism on the
  generic fibre of the ruled surfaces.
\end{itemize}
\end{proposition}
\begin{proof} Let us first prove that
  $S$ and $S'$ are ruled
  surfaces. We already noticed that nef
  terminal models are unique. Hence we
  may assume that $(S,C)$ is not $(\P^2,\overline{C})$. 
  Observe that $\Ka(\P^2,\alpha
  l)=\kappa(\P^2,\alpha l)<0$ for any $\alpha\leq 1$. 
Note that by our construction, when $S$
is a ruled surface and $(S,2/\alpha \overline{C})$
is a model in the list we have
$2/\alpha\leq 1$ and 
\begin{equation}
(K_S+\frac2\alpha\overline{C})\cdot
  f=0,
\label{eq:zero}
\end{equation}
where $f$ is a ruling of $S$. 
Hence we  have
$\overline{\kappa}(S,2/\alpha\overline{C})\geq
0$. This
  shows that $(\P^2,l)$ is not
  birational to other pairs. 

Assume that $\alpha'\geq \alpha$. By
construction
$(S,2/\alpha C)$ has canonical
singularities. Hence $(S,2/\alpha' C)$
has canonical singularities and we have
\begin{equation}
\label{eq:koda}\kappa(S,2/\alpha'
C)=\kappa(S',2/\alpha' C')\geq 0.
\end{equation}
 On
the other hand,  as already observed in
equation (\ref{eq:zero}), every ruled
surface, in our list, satisfies 
$$(K_S+2/\alpha C)\cdot f=0, $$
with $f$ a ruling of $S$. 
This, together with equation (\ref{eq:koda}) yields
$\alpha=\alpha'$.
Assume that $\kappa(S,2/\alpha C)=1$, then
$S\iso\F_a$ and $S'\iso\F_{a'}$. The fibre
structure on $\F_a$ is, again by
equation (\ref{eq:zero}), the log-Iitaka
fibration of the pair. Therefore it is
preserved by any birational map of the pair.
\end{proof}

\begin{remark} Note that
  Proposition \ref{pro:unique} (ii)
  is false for pairs with zero
  log-Kodaira dimension. Consider a
  sestic $C$ with at least four ordinary double
  points. Then blowing up the fourth node
  realizes a model but one can
  equivalently apply a standard Cremona
  transformation on the first three nodes and then
  blow up the fourth.
\end{remark}

\begin{remark}\label{rem:elm}
Let $D\sim \alpha C^a_0+\beta f$ be an irreducible curve in $\F_a$. Let
$elm_p:\F_a\map\F_{a\pm1}$ be the elementary transformation based
on a point $p$ of multiplicity $m$ for the curve $D$. Let
$D'=elm_{p*}(D)\sim a C^{a\pm1}_0+b' f$ the strict transform. Then
according to the position of $p$ with respect to $C_0$ we have
$$\beta'=\beta-m\ \ \mbox{if $p\not\in
  C^a_0$,}\ \ \beta'=\beta+(\alpha-m)\ \  \mbox{if
  $p\in C^a_0$} $$
\end{remark}

\begin{definition} An irreducible and
  reduced curve $C\subset\P^2$ is a
  minimal degree curve if it is not Cremona
  equivalent to any curve of lower
  degree.
\end{definition}

We are ready to put the first brick.

\begin{lemma} Let $C\subset\P^2$ be a
  curve of degree $d>1$, with
  multiplicity set $\{m_i\}$. Assume
  the
  following
\begin{itemize}
\item[a)]  $m_1\geq d/3$
\item[b)]  $C$ is a minimal degree curve
\item[c)] $(S,\overline{C})$ is a model
  of $(\P^2, C)$ in the list of
  Proposition \ref{pro:model}, with
  $S\neq\P^2$ and 
  $\overline{C}\sim\alpha C_0+\beta f$
\item[d)] $\kappa((S,2/\alpha\overline{C})=1$
\end{itemize}
 Then
  $\alpha=d-m_1$,
  $S\iso\F_a$ and the general fibre of $\F_a$ is the
  strict transform of a
  line through $p_1$.
\label{lem:minimal}
\end{lemma}
\begin{remark} Let us stress that the pair
  $(\P^2,3/dC)$, as
  in the statement of Lemma
  \ref{lem:minimal}, has at most two places
  of non canonical singularities on
  $\P^2$. If moreover the places are two
  then all other singularities are
  terminal, again the assumption on
  Kodaira dimension is crucial. We would
  like to thank Alberto Calabri and Ciro
  Ciliberto for pointing this out to us
  while we where cruising in wrong directions.
\end{remark}

\begin{proof} Let $\nu:\F_1\to\P^2$ be
  the blow up of $p_1$. Then
  $$\nu^{-1}_*(C)=\overline{C}\sim
  (d-m_1)C^1_0+d f$$
If either $m_2< \frac{d-m_1}2$  or $m_2>\frac{d-m_1}2$ and
 $m_3<\frac{d-m_1}2$ we produce a
terminal model, therefore unique, and
the claim is clear.

 If
$m_2> \frac{d-m_1}2$ and
 $m_3\geq\frac{d-m_1}2$ then assumption
 b) forces all non terminal points
to lie on the exceptional section. Let  $(\F_a, \tilde{C})$ be a
standard model of $(\P^2,C)$. Let
$\tilde{C}\sim \alpha C^a_0+\tilde{\beta} f$, then
we have.

\begin{claim}The general fibre of
  $\F_a$ is a line through $p_1$ and $\alpha=d-m_1$.
\end{claim}
\begin{proof}[Proof of the claim] The
  standard model $(\F_a,\tilde{C})$ is
  obtained via a chain of elementary
  transformations, starting from the
  blow up of a point of maximal
  multiplicity of $C$. In particular the
  general fibre, after this blow up, is
  a line passing through $p_1$. The only
  ruled surface with two fibre
  structures is $\F_0$. Hence to
  change the general fibre the chain of
  elementary transformations has to go
  back to $\F_0$, and henceforth to
  $\F_1$. Let $\psi:\F_a\rat\F_1$ such a
  map, and
  $$\overline{C}:=\psi_*(\overline{C})\sim (d-m_1)C^1_0+ \overline{d} f$$
be the strict transform of $C$ on $\F_1$. Then by Remark
\ref{rem:elm} we have
$$\overline{d}=d+\sum_1^k((d-m_1)-\mu^0_j)-\sum_1^k
  \mu_h $$
where $\mu^0_j$ are the multiplicities of
points successively lying on the
exceptional section and $\mu_h$ are
the multiplicities of points lying
outside of it. Note that by construction
both $\mu_h$ and   $\mu_j^0$ are bounded
from below by $\frac{d-m_1}2$.
Furthermore by
hypothesis we
have $\mu_1^0=m_2>\frac{d-m_1}2$.
Therefore the existence of such a chain
produces a plane curve of degree $d'<d$ Cremona
equivalent to $C$.
\end{proof}
The claim and  Proposition
\ref{pro:unique} allow to conclude
in this case.

\noindent Assume now that $m_2=\frac{d-m_1}{2}$.
 Let $(\F_a,\tilde{C})$ be a standard
model. Then we
only blow up points on the exceptional
section. In particular no chain of
elementary transformation can land on
$\F_0$. This is enough to conclude again
by Proposition \ref{pro:unique}.
\end{proof}


\begin{definition} Let
  $\Sigma\subset\F_a$ be a (possibly reducible)
  section, $f_1$ and $f_2$ two
  fibres and
  $\tilde\Sigma\sim\Sigma+f$ an
  irreducible section. The linear
  system
$$\Lambda_{\Sigma}=\{\Sigma+f_1,\Sigma+f_2,\tilde\Sigma\}$$
is called a planar linear system
associated to $\Sigma$. Let
$\f:\F_a\rat\P^2$ be the map associated
to $\Lambda_\Sigma$ and $D\subset\F_a$
an irreducible curve. Then $(\P^2,\f_*(D))$
is called the plane model of $(\F_a,D)$
associated to $\Lambda_\Sigma$.
\end{definition}

\begin{remark}\label{rem:planar}
Every rational
  map $\F_a\rat\P^2$ sending the general
  fibre of $\F_a$ to a line is given by
  some planar linear system.
\end{remark}
 As a
  foreword to the next Proposition we want
  to spend few lines on these linear
  systems.
Let
$\Lambda_\Sigma=\{\Sigma+f_1,\Sigma+f_2,\tilde\Sigma\}$
be a planar linear
system on $\F_a$, and $\f:\F_a\rat \P^2$
the rational map associated. The base
locus of $\Lambda_\Sigma$ is given by
$$\Bsl\Lambda_{\Sigma}=\Sigma\cap\tilde{\Sigma}$$
and the morphism $\f$ can be factored as
follows. Let $\theta:\F_a\rat\F_1$ be
the chain of
elementary transformations centered
along $\Bsl\Lambda_\Sigma$. Then we have
$\theta_*(\Sigma)=C_0^1\subset\F_1$. Let
$\nu:\F_1\to\P^2$ be the contraction of
the exceptional section $C^1_0$. In this
notation we have $\f=\nu\circ\theta$. In
particular given an irreducible curve
$D\subset\F_a$ with $D\sim \alpha
C^a_0+\beta f$ then
$$\theta_*(D)\cdot C^1_0=D\cdot\Sigma-\sum \mu_i$$
where the $\mu_i$ are the multiplicities
of $D$ along $\Bsl\Lambda_\Sigma$. Hence
the curve $\f_*(D)\subset\P^2$ has
degree
\begin{equation}
\deg\f_*(D)=D\cdot\Sigma-\sum
\mu_i+\alpha.
\label{eq:degree}
\end{equation}

Our next aim is to single out special
planar systems. Consider a standard model $(\F_a,
\overline{C})$. Let $\nu:\F_a\rat\F_b$ be a chain of
elementary transformations that resolve
the singularities of $\overline{C}$
along the exceptional section.
Let $\hat{C}:=\nu_*(\overline{C})\sim
\alpha C^b_0+\beta' f$
be the strict transform. Note that
$(\F_b,\frac2\alpha\hat{C})$ has not
canonical singularities. Let $x_1\not\in
C^b_0$ be a point and
$elm_{x_1}:\F_b\rat \F_{b-1}$ the
elementary transformation centered at
$x_1$. Recursively choose a point $x_i\in \F_{b-i-1}\setminus C^{b-i-1}_0$ and
perform the elementary transformation
centered at $x_i$. Then after $b-1$
steps we have
$$\psi_{x_1,\ldots
  x_{b-1}}:=\epsilon\circ elm_{x_{b-1}}\circ\cdots\circ
elm_{x_1}:\F_b\rat \P^2$$
where $\epsilon$ is the blow down of the
$(-1)$-curve on $\F_1$. In particular
$(\psi_{x_1,\ldots
  x_{b-1}})_*(\hat{C})$ is a curve of
degree
\begin{equation}
\label{eq:section}
\hat{C}\cdot(C^b_0+(b-1)f)-\sum_1^{b-1}\mu_i+\alpha=\beta'
-\sum_1^{b-1}\mu_i
\end{equation}

where $\mu_i=\mult_{x_i}\hat{C}$.
\begin{definition}\label{def:mpm} Let
$(\F_a, \overline{C})$ be a model of Proposition
  \ref{pro:model}. Let $\nu:\F_a\rat\F_b$ be a chain of
elementary transformations that resolve
the singularities of $\overline{C}$
along $C_0$.
Let $\hat{C}:=\nu_*(\overline{C})$
be the strict transform.
Let $\{x_1,\ldots, x_{b-1}\}$ be a set of
points  as in the above
construction. We say that $\{x_1,\ldots,
x_{b-1}\}$ is minimal for $\hat{C}$ if $\sum\mult_{x_i}\hat{C}$ is
maximal.
 Consider the
planar linear system on $\F_b$
$$\Lambda_{C^b_0}:=\{(C^b_0+\sum_1^{|b-1|}f_i)+f_0,
(C^b_0+\sum_1^{|b-1|}f_i)+f_b,\tilde{\Sigma}\}$$
where $\{\tilde\Sigma\cap
f_i\}_{i=1}^{b-1}$, is a minimal set for
$\tilde{C}$.
Then the plane model associated to
$\Lambda_{C^b_0}$ will  be called a
  minimal plane model for
  $(\F_a,\overline{C})$.
\end{definition}
\begin{remark} Unfortunately the
  definition of minimal plane models
  requires a partial resolution of
  singularities. Note that, in general,
  this cannot be
  avoided. To convince yourself it is
  enough to  consider a
  standard model with a unique, nested,
  singularity along the exceptional
  section. 
The hidden
  contribution coming from the
  singularity could drop the degree much
  more than the smooth points
  outside, cfr.
  Example \ref{ex:bad}
\end{remark} 
 With this in mind we are ready for the
following\begin{theorem}
\label{th:main}Let $C\subset\P^2$ be a curve, with
multiplicity set $\{m_i\}$. The curve
$C$ is a minimal degree curve if and
only if either $m_1< d/3$ or it
is a minimal plane model or it is a line.
\end{theorem}

\begin{proof}
 Let $(\P^2,C_{\min})$ be a
  minimal degree curve. Assume that
  $C_{\min}$ is not a line and  $m_1\geq
  d/3$.
 Then a standard model is a ruled
  surface. Let
  $(\F_a,\overline{C})$ be a standard
  model for $(\P^2,C_{\min})$. Let
  $\nu:\F_1\to\P^2$ be the blow up of a
  point
  $p\in C_{\min}$ of maximal
  multiplicity. Let
  $\tilde{C}=\nu^{-1}_*(C_{\min})$ be the
  strict transform, and $E\subset\F_1$
  the exceptional section.

Let $\overline{C}\sim\alpha C^a_0+\beta
f$, and assume $\kappa(\F_a,2/\alpha \overline{C})=1$.
Then,
  by Lemma  \ref{lem:minimal} and its
  proof, there is
  a sequence of elementary
  transformations $\Phi:\F_1\rat \F_a$, leading
  $(\F_1,\tilde{C})$ to
  $(\F_a,\overline{C})$, with
  $\Phi_*(E)=C_0^a$. Let
  $\chi:\F_a\rat \F_b$ be a resolution of
  the singularities of $\overline{C}$ along
  $\Phi_*(E)$, with $\hat{C}=\chi_*(\overline{C})\sim \alpha C^b_0+\beta' f$.
Then the linear system taking
  $(\F_b,\hat{C})$ to $(\P^2,C_{\min})$ is a
  planar linear system associated to
  $\Lambda_{C_0^b}$, keep in mind
  Definition \ref{def:mpm}. Hence
  $C_{\min}$ is a minimal plane model.

Assume that $\kappa(\F_a,2/\alpha
\overline{C})=0$ then  $a\leq 2$.

\begin{claim} A standard model of
  $(\P^2,C_{\min})$ is obtained with at
  most one elementary transformation.
\label{cl:0}
\end{claim}
\begin{proof}[Proof of the Claim] Let
  $\nu:\F_1\to\P^2$ be the blow up of
  $p_1\in C_{\min}\subset\P^2$ and
  assume that
  $(\F_1,\nu^{-1}_*(C_{\min}))$ is not
  the model. If we have
  either $p_2\not\in C^1_0$ or $p_3\not\in
  C^1_0$. Then minimal degree forces $\F_a$
  to be either $\F_0$ or $\F_2$ with
  $m_3\leq (d-m_1)/2$, and
  the claim is clear.

Assume $p_2,p_3\in C^1_0$, and $m_3\geq
(d-m_1)/2$. The requirement on Kodaira
dimension forces every model to have
$a\leq 2$. Therefore, after maybe
the first,  any elementary
transformation centered on the
exceptional section has to be balanced
by an elementary transformation
centered outside the exceptional
section. If this chain of elementary
transformations take us back onto $\F_2$
with a curve $\tilde{C}\subset\F_2$  then
$$\tilde{C}\sim(d-m_1)C^2_0+(d+\sum_1^s((d-m_1)-\mu^0_j)-\sum_1^{s-1}
\mu_h) f $$
where $\mu_j^0$ and $\mu_h$ are the
multiplicity of points on $C^2_0$,
respectively, outside $C^2_0$.
The standard model construction yields
$\mu^0_j\geq (d-m_1)/2$
and $\mu_h>(d-m_1)/2$. If
$(\F_2,2/\alpha\tilde{C})$ is not
terminal the projection
from a non terminal point produces a
curve of degree
$$\tilde{d}\leq d+\sum_1^s((d-m_1)-\mu^0_j)-\sum_1^{s-1}
\mu_h-(d-m_1)/2<d $$

If $(\F_2,2/\alpha \tilde{C})$ is terminal then it is the standard
model and $\tilde{C}\cdot C^2_0=0$. In particular $\mu_h=d-m_1$, for
any $h=1,\ldots,s-1$, hence a general projection of $\tilde{C}$
onto $\P^2$ has degree
$$\tilde{d}=d+\sum_1^s((d-m_1)-\mu^0_j)-(s-1)(d-m_1)<d. $$
\end{proof}

\noindent The claim is enough to conclude that
$C_{\min}$ is a minimal plane model of
$(\F_a,\overline{C})$.

We need to prove that the curves
described are of minimal degree.
If $(\P^2,C)$ satisfies
$m_1<d/3$, then  it is a minimal degree curve
by Lemma \ref{lem:NF}, a line is clearly
of minimal degree. 
Assume that $(\P^2,C)$ has $m_1\geq
d/3$, it is not a line and it is a minimal
plane model. Let $(\F_a,\overline{C})$
be a standard model, with
$\overline{C}\sim\alpha C^a_0+\beta f$.

We first prove that minimal plane models
are minimal degree curves between all
planar models of $(\F_a,\overline{C})$.

\begin{claim} Minimal plane models of
  $(\F_a,\overline{C})$ have minimal
  degree among all planar models of $(\F_a,\overline{C})$.
\end{claim}
\begin{proof}[Proof of the Claim]
Let
$$\Lambda_{\Sigma}=\{\Sigma+f_1,\Sigma+f_2,
  \tilde{\Sigma}\}, $$
be a planar linear system, with $\Sigma\sim C^b_0+\gamma f $.
Then we have
$$\Sigma\cdot\tilde\Sigma=2\gamma+1-b,$$
and, by equation (\ref{eq:degree}), the plane model associated to
$\Lambda_{\Sigma}$ has degree
\begin{equation}\delta_{\Sigma}=\Sigma\cdot
\hat{C}-\sum_1^{2\gamma+1-b} \mu_i+\alpha=
\alpha(\gamma-b+1)+\beta'-\sum_1^{2\gamma+1-b}
\mu_i
\label{eq:deg}
\end{equation}
where the $\mu_i$ are multiplicities of
$\hat{C}$ at the centres of elementary
transformations.

To prove the claim we have to
check that the degree of the plane curve produced
by the  planar linear system associated to the section
$\Sigma_{C^b_0}=C^b_0+(b-1)f$ is at most
$\delta_\Sigma$.
Hence, keeping in mind equations
(\ref{eq:section}), and (\ref{eq:deg}), we have  to check that
$$-\sum_1^{b-1}\mu_i^0\leq\alpha(\gamma-b+1)-\sum_1^{2\gamma+1-b}
\mu_i $$ where $\mu_i^0$ are the multiplicities associated to the
planar linear system
$\Lambda_{\Sigma_{C^b_0}}$, and $\mu_i$
are the multiplicities associated to
the planar linear system $\Lambda_\Sigma$. There is
a unique curve with negative self
intersection on $\F_1$. Hence  at least
$b-1$ points of $\Sigma\cdot \tilde\Sigma$ are outside the
exceptional section. Moreover $(\F_a,2/\alpha \overline{C})$ is
a standard model therefore there are at most $b-a$ points of
multiplicity strictly greater than
$\alpha/2$, and this points are 
outside the exceptional curve. That is we
may  assume, after reordering the indexes,
$$\sum_1^{b-1}\mu_i^0\geq \sum_1^{b-1}\mu_i$$
Then we have to check that
$$\alpha(\gamma-b+1)-\sum_{b}^{2\gamma+1-b}
\mu_i \geq 0$$ Moreover, we have that
$$\mu_i\leq \frac\alpha2\ \ \ \mbox{for
  $i> b-a$}$$
then
$$\alpha(\gamma-b+1)-\sum_{b}^{2\gamma+1-b}
\mu_i \geq\alpha(\gamma-b+1)-2(\gamma-b+1)\alpha/2=0. $$

\end{proof}

Next we prove that it is enough to
consider planar models of $(\F_a,\overline{C})$.

\begin{claim} We may assume that  any minimal degree plane
curve birational to
$(\F_a,\overline{C})$ is obtained via a
planar model of $(\F_a,\overline{C})$.
\end{claim}
\begin{proof}[Proof of the Claim]
If
$\kappa(\F_a,2/\alpha\overline{C})=1$,
we know by Lemma \ref{lem:minimal}, and
Remark \ref{rem:planar}, that any 
minimal degree plane curve birational to
$(\F_a,\overline{C})$ is obtained via a
planar linear system. 

Assume that
 $\kappa(\F_a,2/\alpha\overline{C})=0$ and
let $(\F_a,\overline{C})$ and
$(\F_{b}, \tilde{C})$ be two birational models in the
list of Proposition \ref{pro:model}. We
may assume, without loss of generality,
 that $a=1$ and $b=2$.
The vanishing of Kodaira dimension
and $\alpha$ uniquely determine the
linear equivalence class of
$\overline{C}$ and $\tilde{C}$, namely
$$\overline{C}\sim \alpha
C^a_0+\frac32\alpha f\ {\rm
  and}\ \tilde{C}\sim \alpha
C^b_0+2\alpha f  $$
Note further that the existence of two
models forces the presence of a canonical
singularity. Therefore the minimal plane
models obtained by the two have equal
degree $\frac32\alpha$. 
\end{proof}

The two claims conclude the proof of the Theorem.
\end{proof}

The main question is: can we describe  minimal
plane model curves (without going
through a partial resolution)?
We do not expect to have a positive
answer in general.
In the positive direction there is a nice result of
Jung, we are able to recover.

\begin{corollary}[\cite{Ju}] Let $C\subset\P^2$ be
  a curve with $m_1+m_2+m_3\leq d$. Then
  $C$ is a minimal degree curve.
\end{corollary}
\begin{proof} The assumption forces
  $m_i\leq (d-m_1)/2$. Therefore a
  minimal plane model is reached by
  performing the inverse of a resolution
  of singularities along the exceptional
  section.
\end{proof}
Unfortunately the opposite direction is not true, even
discarding the trivial example of lines.

\begin{example} Let $C\subset\P^2$ be a
  curve of degree $7$ with a point of
  multiplicity $4$ and two infinitely
  near double points. Then it is easily
  seen that $C$ is a minimal curve.
\end{example}

The difficulty in predicting minimality
for plane curves can be seen in the
following example.

\begin{example} Let $D_i\subset\F_3$ be
  irreducible and reduced curves with
  $D_i\sim 3C^3_0+11f$. Assume that $D_1$
  has a unique ordinary double point,
  say $p_1\in C^3_0$ and $D_2$ has a
  unique ordinary double point, say
  $p_2\in\F_3\setminus C^3_0$. The pair $(\F_3,D_i)$
  is  birational to a
  pair $(\P^2,C_i)$. Where $C_i$ is a degree
  $9$ curve with a six-tuple point and $3$
  infinitely near double points. The
  main difference, that we can only
  divine on the resolution, is that $C_1$
  is a minimal degree curve, while $C_2$
  is birational to a curve of degree $8$
  with a quintuple point and an
  infinitely near double point.
\label{ex:bad}
\end{example}
It is easy to produce similar examples
of arbitrarily nested singularities.
\section{Rational surfaces}

It is natural to ask when a rational hypersurface is Cremona
equivalent to a hyperplane. 
The unique case in which an answer to
the above question is
known is that of rational curves. Coolidge, \cite{Co}, first
suggested that this should be the case if and only if
\hbox{$\Ka(\P^2,C)<0$}.  
Kumar--Murthy, \cite{KM} see also
\cite[Proposition 12]{Ii}, were able to
prove Coolidge statement and refine it a
little bit. Approaching the question
via LMMP we found an alternative proof.

\begin{proposition} A rational plane curve
  $C$ is Cremona equivalent to a line if
  and only if $\Ka(\P^2,1/2 C)<0$.
\label{pro:2mmp}
\end{proposition}
\begin{proof} We have only to prove that
  if $\Ka(\P^2,1/2 C)<0$ then $C$ is
  Cremona equi\-va\-lent to a line. Let
  $g:S\to\P^2$ be a minimal resolution
  of singularities, with
  $C_S=g^{-1}_*(C)$.  Then by hypothesis
  we have
  $\kappa(S,1/2 C_S)<0$. Let us start a
  LMMP program for the pair $(S,1/2
  C_S)$. Then by construction   every
  curve contracted during the LMMP satisfies
$$(K+1/2 \overline{C})\cdot E<0$$
where $K$ is the canonical class and $\overline{C}$ the strict
transform of the curve $C$. This means
$C\cdot E<2$. In particular the output
of a LMMP is a pair $(\overline{S},1/2\overline{C})$
with $\overline{C}$ smooth and $K_S+1/2\overline{C}$ negative on
infinitely many curves. If $\overline{S}$ is $\P^2$ the claim is
clear. Otherwise $\overline{S}$ is a
rational ruled surface, 
$\overline{C}$ is a smooth rational
curve and 
\begin{equation}
(K_S+1/2 \overline{C})\cdot f<0,
\label{eq:meno}
\end{equation}
where $f$ is a fiber in a ruling of $S$.
If $\overline{C}$ is a section of a
ruling, or
a fiber then the claim is clear. Assume
that this is not the case and let
$S\iso\F_a$, and $\overline{C}\sim
\alpha C_0^a+\beta f$, for some
$\alpha\geq 2$. Then by equation
(\ref{eq:meno}) we have $2\leq
\alpha\leq 3$, and genus formula gives
\begin{eqnarray*}-2=((-2+\alpha)C_0+(-2-a+\beta)f)\cdot(\alpha
C_0+\beta
f)=\\
=\alpha(a(\alpha-2)+(\beta-a-2))+\beta(\alpha-2)
\end{eqnarray*}
The assumption that
$\overline{C}$ is irreducible and it is not a section shows
that the only possibility is
$a=1$ and $\alpha=\beta=2$. That is
 $(S,\overline{C})$ is equivalent to a
 plane conic and therefore to a line. 
\end{proof}

\begin{remark} The proof shows that the
  numerical hypothesis on log-Kodaira
  dimension allows to choose a Mori
  fiber Space structure
  particularly ``nice'' for the curve $\overline{
C}$.
\end{remark}
The case of plane curves has been studied
from many different points of view.
Various other necessary and sufficient
conditions are known, see \cite{BB} for
a nice survey. Probably the
most tempting conjecture is
Nagata--Coolidge's prediction that every
cuspidal rational curve is Cremona
equivalent to a line, \cite{Co} \cite{Na}.
For this we do not see any
translation into LMMP dictionary.

It is quite natural to ask for generalisations in higher
dimensions. The main difficulty is the poor knowledge of Cremona
group starting from $\P^3$.

The case of surfaces is already
quite mysterious. It is easy to show
that quadrics and rational cubics are Cremona
equivalent to a plane, see 
Case \ref{case:3} below. Rational
quartics with either 3-ple or 4-uple
points are again easily seen to be
Cremona equivalent to planes, the latter
are cones over rational curves Cremona
equivalent to lines. 

 It has been expected that
Noether quartic should be the first
example of a rational surface
not Cremona equi\-va\-lent to a plane, but this is not the case.

\begin{example} Let $S\subset\P^3$ be
  the Noether quartic. That is a quartic
  with a unique double point of local
  analytic type
  $O\in(x^2+f_3+g_4=0)\subset\C^3(x,y,z)$. We
  may assume that the equation of $S$ is
$$(x_0^2x_3^2+f_3x_3+g_4=0)\subset\P^3$$
with $p\equiv[0,0,0,1]\in S$ the unique
singular point.
Let
   $\epsilon:Y\to\P^3$ be a weighted
blow up of $p$,
  with weights $(2,1,1)$ on the
  coordinates $(x_0,x_1,x_2)$, and
  exceptional divisor $E$. Then we have
  $\epsilon^*(x_0=0)=H+2E$ and
  $\epsilon_{|H}:H\to (x_0=0)$ is an
  ordinary blow up. Let
  $l\subset(x_0=0)$ be a line through
  $p$, and $l_Y$ its strict
  transform. Then we have $H\cdot l_Y=-1$ and
  $K_Y\cdot l_Y=-1$. This shows that we can
  contract $H$ to a curve. Let $\mu:Y\to
  V$ be the contraction of $H$. To
  understand $V$ let us first observe
  that $\rk Pic(V)=1$. Then consider
  a general line $p\in
  r\subset\P^3$ and $r_Y$, respectively
  $r_V$ its strict transforms. By construction 
\begin{equation}\label{eq:noether}K_V\cdot r_V=K_Y\cdot
  r_Y=(\epsilon^*(K_{\P^3})+3E)\cdot
  r_Y=-4+\frac32=-\frac52
\end{equation}
This shows that $V$ is a $\Q$-Fano
3-fold of Fano index at least $5/2$ and
with a singular point of type
$1/2(1,1,-1)$. It
is well known, see for instance \cite{CF}, that 
$V$ is the cone over the Veronese
surface. 
Further note, for future reference, that equation
(\ref{eq:noether}) also shows that
$\Ka(\P^3,t S)<0$ for any $t<5/4$.

For what follows it is
important to understand this birational
modification in terms of linear systems
of $\P^3$. Consider the linear
  system $\Lambda\subset|\O(2)|$, 
of quadrics with multiplicity $2$ on the
  valuation $E$. That is quadrics, say $Q$,
  through $p$ with
  ${\mathbb T}_pQ\supseteq(x_0=0)$. Then
  we have
  $\dim\Lambda=6$ and
  $V:=\f_\Lambda(\P^3)\subset\P^6$ is the
  cone over the Veronese surface in
  $\P^5$. Moreover the strict transform
  of the surface $S$
$$S_V:=\f_\Lambda(S)_*=(\mu\circ\epsilon^{-1})_*(S)\in|\O_V(2)|,$$
  is a smooth surface in $V$.

Our next aim is to project the pair
$V$ birationally onto a cubic
3-fold $T\subset\P^4$ in such a way that
$S_V$ is sent to an hyperplane section.
To do this we have to produce a
birational linear
system, say ${\mathcal C}$, of dimension
4 and containing $S_V$.

Let  $A_1$ be a general hyperplane section
  of $V$ and consider $\Gamma:=A_1\cap
  S_V$, a smooth irreducible
  curve. Fix $C\subset S_V$ a general
  smooth conic, and let
  $$\B=\{B_0,B_1,B_2,B_3\}=|\O_V(1)\otimes\I_C|,$$
be the linear system of hyperplanes
containing $C$.

The desired linear system is then
$${\mathcal
  C}=\{S_V,A_1+B_0,A_1+B_1,A_1+B_2,A_1+B_3\}\subset|\O_V(2)|$$ 
Let $\f=\f_{\mathcal C}:V\rat \P^4$ the
rational map associated. 

We have to prove that $\f$ is birational
onto the image and compute the degree of $\f(V)$.
To do this we have to understand the
intersection of elements in ${\mathcal
  C}$.

By construction 
$$
S_V\cdot(A_1+B_i)=\Gamma+C+R_i $$ for some residual curve
$R_i\subset B_i$ and
$$B_{j|B_i}=C+D_{ji}, $$
for some residual curve $D_{ji}$.
 Let $v:\P^2\to
B_i$ be the 2-Veronese embedding. 
In this notations $v^*(C)$ is a line,
$v^*(S_{V|B_i})\sim\O(4)$ and $v^*(B_{j_|B_i})\sim\O(2)$.
This yields the following:
\begin{itemize}
\item[(a)]  $R_i$  is the image of a cubic curve
via the Veronese map $v$;
\item[(b)] $v^*(D_{ji})\sim\O(1)$.
\end{itemize}
This gives  $(D_{ji}\cdot
D_{hi})_{B_i}=1$, and $(R_i\cdot
D_{ji})_{B_i}=3$. Hence  $\f$ is
birational and $\deg\f(V)=3$. 
This means that $\f(V)$ is a rational,
hence singular, cubic hypersurface and
$\f(S_V)$ is an hyperplane section. Then
the projection from a singular point of
$\f(V)$ realizes $\f(S_V)$ as a rational
surface of degree at most three in $\P^3$. This is enough to conclude.
\label{ex:noether}
\end{example}

 This
example and the proof of Proposition \ref{pro:2mmp} suggest that a
result similar to Proposition \ref{pro:2mmp} is at hand also for
$\P^3$. To get it we have to slightly modify the $\sharp$-MMP
developed in \cite{Me} for linear systems on uniruled 3-folds.
First we recall the notion of
(effective) threshold.
\begin{definition}
Let $(T,H)$ be a $\Q$-factorial uniruled
   3-fold and $H$ an irreducible and
   reduced Weil divisor
  on $T$. Let
  $$\rho=\rho_{H}=\rho(T,H)=:\mbox{  \rm sup   }\{m\in \Q|H+mK_T
\mbox{
  \rm is an effective $\Q$-divisor  }\}\geq 0,$$
 be the (effective) threshold of the pair
 $(T,H)$.
\end{definition}
\begin{remark} The threshold is not a
birational invariant of the pair and it
is not preserved by blowing up.
Consider
 a quadric cone $Q\subset\P^3$ and let
$Y\to \P^3$ be the blow up of the vertex
then $\rho(Y,Q_Y)=0$, while
$\rho(\P^3,Q)>0$. 
\end{remark}
Keeping in mind the surface case we want
to use the threshold to single out a
Mori fiber Space
after a LMMP and then use this
``simplified'' model to prove the
Cremona equivalence to a
plane. To do this we need to start from
a pair with mild singularities, say $(Y,S_Y)$, and
run a LMMP in a way to get a nice MfS
structure with respect to the strict
transform of $S$.

\begin{remark} Let us stress a point
  where the 3-fold and surface cases are
  quite different. Assume that
  $\rho(Y,S_Y)=0$. In this case $S_Y$ is
  birationally trivial along fibers of the
  MfS associated. While for surfaces
  this is an easy case
 for 3-folds it is extremely hard
  to extract enough
  information from the Mori fiber space
  associated. 
The main
  difficulty comes from rational conic
  bundles, that are not scrolls. The knowledge of these 3-folds is very poor.
Even assuming the standard conjectures, like Cantor or Iskovskikh rationality criteria, it is very difficult to guess any
  kind of birational embedding of  $S_T$
  in $\P^3$. 
\end{remark}
 
The remark suggests to introduce the
following variants.

\begin{definition}
Let $(Y,S_Y)$ be a 3-fold pair birational
to a 3-fold pair $(T,S)$. We say that
$(Y,S_Y)$ is a good birational model if $Y$ has
terminal $\Q$-factorial singularities
and $S_Y$ is a smooth Cartier divisor. 
The sup-threshold
of the pair $(T,S)$ is
$$\overline{\rho}(T,S):=\sup\{\rho(Y,S_Y)\},$$
where the sup  is taken on good
birational models.
\end{definition}

\begin{remark} It is clear that any pair
  $(\P^3,S)$ Cremona equivalent to a
  plane satisfies
  $\overline{\rho}(\P^3,S)>0$. The pair
  $(\P^3,H)$, where $H$ is a plane, is a
  good model with positive
  threshold.\label{rem:pos}

Considering birationally super-rigid
MfS's one can produce examples of
pairs, say $(T,S)$, with
$\overline{\rho}(T,S)=0$. It is not
clear to us if such examples can exist
also on 3-folds with bigger pliability,
\cite{CM}.  
\end{remark}

The main result we need in this contest
is the following modification of
\cite[Theorem 5.3]{Me}.

\begin{theorem}\label{th:3mmp}
 Let $(\P^3,S)$ be a pair with $S$ a
 rational surface. Assume that for some
 good model $(T,S_T)$ the divisor
 $K_T+S_T$ is not pseudoeffective and 
 $\overline{\rho}(\P^3,S)>0$. 
Then $(\P^3,S)$ has a good birational
model in the following list:
  \begin{itemize}

  \item[i)] a rational $\Q$-Fano $3$-fold $T$ of index $1/\rho>1$,
with $K_{T}\sim-1/\rho S_T$, of the
following type:
\begin{itemize}
\item[a)] $(\Proj(1,1,2,3),\O(6))$
\item[b)] $(X_6\subset \Proj(1,1,2,3,a), X_6\cap \{x_4=0\})$, with $3\leq a
\leq 5$
\item[c)] $(\Proj(1,1,1,2),\O(2k))$,
  with $k=1,2$,
\item[d)] $(X_4\subset \Proj(1,1,1,2,a), X_4\cap \{x_4=0\})$, with
$2\leq a \leq 3$ \item[e)] $(\Proj^3,\O(a))$, with $a\leq 3$,
\item[f)] $(X_3\subset \Proj(1,1,1,1,2),X_3\cap \{x_4=0\})$,
\item[g)] $(\Q^3,\O(b))$, with $b\leq 2$, $X_{2,2}\subset
\Proj^5,\O(1))$, a linear section of the Grassmann variety
parametrising
  lines in $\Proj^4$, embedded in $\Proj^9$ by Plucker
coordinates;
\end{itemize}
\item[ii)] a bundle over  $\P^1$
  with generic fibre $(F,S_{T|F})
\iso(\Proj^2,\O(2))$ and with at most finitely many fibres
$(G,S_{T|G})\iso (\s_4,\O(1))$, where $\s_4$ is the cone over the
normal quartic curve and the vertex sits over an hyper-quotient
singularity of type $1/2(1,-1,1)$ with $f=xy-z^2+t^k$, for $k\geq
1$, \cite{YPG};
 \item[iii)]
  a quadric bundle over $\P^1$ with at most
$cA_1$ singularities of type $f=x^2+y^2+z^2+t^k$, for $k\geq 2$,
and $S_{T|F}\sim \O(1)$;
 \item[iv)]$(\Proj(E),\O(1))$ where $E$ is
a $\rk 3$  vector bundle over $\P^1$;
  \item[v)] $(\Proj(E),\O(1))$ where $E$ is a $\rk 2$  vector
bundle over a rational surface $W$;
  \end{itemize}
\end{theorem}

\begin{proof} The hypothesis and Lemma
  \ref{lem:psef} yield that $K_Y+S_Y$
is not pseudoeffective for any good
birational model. That is
$\rho(Y,S_Y)<1$. 

Fix a good birational model $(Y,S_Y)$
with $\rho(Y,S_Y)>0$. We mimic the proof
  of \cite[Theorem 5.3]{Me}.  The main difference
  is that we are assuming that the
  3-fold is rational
  and $S$ is not a big linear system but
  a fixed divisor.

Let us start a LMMP for the pair
$(Y,(1-\epsilon)S_Y)$, for
$0<\epsilon\ll1$.
Let $[C]$ be the class of an extremal
ray encountered along the program. If
$C\cdot S_Y\geq 0$ then this is a usual
step of $\#$-MMP and the birational
transformation associated produce a new
good birational model $(Y^1,S^1)$, see \cite[Proposition 3.6]{Me}. Assume that
 $S_Y\cdot C<0$ and let $\f:Y\rat Y'$
be the birational modification
associated. Then all curves numerical
proportional to $C$ are contained in $S_Y$.
By construction we have
$$0<\rho(Y,S_Y)\leq \rho(Y',\f_*S_Y),$$
the 3-fold $Y$ is rational. Hence
$S_Y$ is not $\f$-exceptional and
$\f_{|S_Y}$ is a well defined morphism.
In particular $\f$ is an isomorphism in
codimension 1. By hypothesis
$(K_Y+S_Y)\cdot C<0$ and the normality
of $S_Y$ yields that $(C\cdot
C)_{S_Y}<0$. This means that the
exceptional locus of $\f_{|S_Y}$ is a
bunch of disjoint $(-1)$-curves on $S_Y$. Let
$C_1$ be an irreducible
component of this exceptional locus, and $\nu:Z\to Y$ the blow up of
$C_1$ with exceptional divisor $E$. By
construction we have
$\nu^*S_Y=S_Z+E$. Let $\Gamma=S_Z\cap
E$, then we have
$$E\cdot\Gamma=(\Gamma\cdot\Gamma)_{S_Z}=-1,\ {\rm
  and}\ 
(\Gamma\cdot \Gamma)_E=S_Z\cdot
\Gamma=S_Y\cdot C-E\cdot\Gamma\leq 0,$$
where for the latter inequality we used
 that $S_Y\cdot C$ is a
negative integer.
This means that $\Gamma$ is either a
ruling of $E$ or the exceptional
section. In other words
$N_{C_1/Y}\iso\O(-1)\oplus\O(-a)$, for
some integer $a>0$. Hence $\f$ is a
generalized Francia's antiflip (or
Atiyah's flop if
$a=1$) and $\f(S_Y)$ is a smooth surface
on the smooth locus of the terminal  $\Q$-factorial
3-fold $Y'$. 
This shows that after finitely many
birational modifications we end up with
a good birational model in the list in
 \cite[Theorem 5.3]{Me} with $T$
rational. In particular for Fanos it is
enough to amend the list in
\cite{Me} using the list of non rational
smooth Fanos in \cite[Theorem
  5]{Isk}. This gives all cases from i)
through v).
\end{proof}

Let us start to prove some Cremona equivalences.

\begin{lemma} Assume that $(T,S_T)$ is
  in cases {\rm ii)--v)} of the above
  list. Then the surface $(T,S_T)$ is
  birational to $(\P^3,H)$, where $H$ is
  a plane.
\label{lem:fiber}
\end{lemma}

\begin{proof} We treat the different
  possibilities separately.
\begin{case}[ii] Let $\eta$ be
the general point of the base $\P^1$. Consider
$F_\eta\iso\P^2_\eta$, the generic fibre of $T$, and $C_\eta$ the
generic fibre of $S_T$, a conic. Then  by Tsen's theorem $C_\eta$
has many rational points. Choose three general points on $C_\eta$
and do a standard Cremona transformation on $F_\eta$ centered at
these points. This sends $F_\eta$ to $\P^2_\eta$ and $C_\eta$ to a
line $l_\eta$. This induces a birational modification on $T$ and
shows that $(T,S_T)$ is birational to a pair  in case iv).
\end{case}

\begin{case}[iii] As in the
previous case let $\eta$ be the generic
point of the base $\P^1$. Then
$F_\eta\iso\Q_\eta\subset\P^3_\eta$ is the generic fibre
of $T$ and $C_\eta$ is again a conic
with many points. The projection from a
point sends $F_\eta$ to $\P^2_\eta$ and
$C_\eta$ to a line $l_\eta$. Therefore,
also in this case the pair $(T,S_T)$
 is birational to a pair in
case iv).
\end{case}

\begin{case}[iv]\label{case:iv} Via
elementary transformations centered in a
bunch of points we can modify $(T,S_T)$
into the pair
$(\P(\O^{\oplus2}\oplus\O(-1)),\O(1))$. This
shows that $(T,S_T)$ is birational to a
surface $S\subset\P^3$ of degree $d$
with a line of multiplicity $d-1$. Let
$l\subset S$ be this line.
To conclude these cases we prove by induction on the
degree $d$ that $S$ is Cremona equivalent to
a plane.  The case $d=2$ is
  immediate. Let $\omega:\P^3\rat\P^3$
  be the $T_{2,3}$ map,  given by the linear
  system
$$|\O(2)\otimes\I_l\otimes_1^3\I_{p_i}|$$
where $p_i$ are general points of $S$. Let $l'$ be the image of
the plane $\langle p_1,p_2,p_3\rangle$ via $\omega$, keep in mind
Construction \ref{con:23}. Then $\omega(S)\subset\P^3$ is a
surface of degree 
$$3d-2(d-1)-3=d-1$$ with the line $l'$ of
multiplicity $2d-(d-1)-3=d-2$. We conclude by induction
hypothesis.
\end{case}

\begin{case}[v] Let $D\subset
T$ be a general very ample divisor, and
$$\Gamma:=S_T\cap D$$
the smooth intersection, keep in mind that $S_T$ is smooth. Let
$\elm_\Gamma:T\rat T_1$ be the elementary transformation centered
on $\Gamma$. Let $S_1:=\elm_\Gamma(S)$, then as we observed in
Construction \ref{con:scroll} $T_1$ is still a scroll over the
rational surface $W$, and $S_1\iso W$.
 Let $\nu:W\rat \P^2$ be
a birational map. Then, as explained in
\cite[5.7.4]{Me}, we can follow this
birational map on the 3-fold. This
produces a birational map $\phi:T_1\rat
T'$ and diagram
$$\xymatrix{
T_1\ar[d]_\pi\ar@{.>}[rr]_{\phi}&&T'\ar[d]^{\pi'}\\
W\ar@{.>}[rr]_{\nu}&&\P^2}
$$
In this
way we end up with a pair $(T',S')$,
birational to $(T,S_T)$, with
$T'$ a scroll over $\P^2$ and $S'\iso
\P^2$ a section. Then
$S'_{|S'}\sim\O(s)$ for some
$s\leq0$. If $s<0$ consider  a
curve $C\subset T'$ with
$\pi^{\prime}_*(C)\sim\O(s)$ and $C\cap
S'=\emptyset$. Then $\elm_C$  modify
$(T',S')$ into
$(\P^2\times\P^1,F)$, with $F=p^*(\O_{\P^1}(1))$. To conclude it
is then enough to project from two general
points of this 3-folds in its Segre
embedding.
\end{case}
\end{proof}

We are ready to characterize surfaces
Cremona equivalent to a plane.

\begin{theorem} Let $S\subset\P^3$ be an
  irreducible and reduced surface.
The following are equivalent:
\begin{itemize}
\item[a)]  $S$ is Cremona equivalent to a plane,
\item[b)]  $\overline{\rho}(T,S)>0$ and there is a good model
  $(T,S_T)$ with $K_T+S_T$ not pseudoeffective.
\end{itemize}
\label{th:3c}
\end{theorem}
\begin{remark} The difficult part of the
  above criteria is the bound on the
  sup-threshold. It could
  be difficult to check it for surfaces $S\subset\P^3$
  such that on a log resolution the
  threshold vanishes.
\end{remark}

\begin{proof} One direction is clear,
  keep in mind remark \ref{rem:pos}.

Assume that the conditions in b) are satisfied. 
Then $(\P^3,S)$ has a good model, say $(T,S_T)$ in the
list of  Theorem
  \ref{th:3mmp}. If $(T,S_T)$ is in
  cases ii)-v) then $S$ is Cremona
  equivalent to a plane by  Lemma
  \ref{lem:fiber}.

Assume that $T$ is a Fano as in case i)
of Theorem \ref{th:3mmp}.
We do not have a general
argument. To conclude we have to produce
for any variety in the list a birational
map to $\P^3$ that sends $S_T$ to a
surface we are able to treat.
Let $(T,S_T)$ be the pair birational to $(\P^3,S)$.
\begin{case}[a, b] Assume first
   $(T,S_T)=(\P(1,1,2,3),\O(6))$.
Let $\epsilon:Y\to\P^3$
  be the weighted blow up of the point
  $[0,0,0,1]$ with weights
  $(1,2,3)$, on the coordinates
  $(x_0,x_1,x_2)$,  and exceptional divisor
  $E$. We choose the coordinate in such a
  way that the plane $H=(x_2=0)$ satisfies
  $$\epsilon^*H=H_Y+3E,$$
that is  $\epsilon_{|H_Y}:H_Y\to H$ is the
  weighted blow up $(1,2)$. In
  particular  on $Y$
  there is a line $l$, the strict
  transform of $(x_1=x_2=0)$, with
  normal bundle
$$N_{l/Y}=\O(-1)\oplus\O(-2)$$
Let $\mu:Y\rat Y'$ be the antiflip of
$l$, then $H'=\mu_*H_Y=\P(1,1,2)$ is a cone with
the vertex over a terminal point of type
$1/2(1,1,-1)$ and we can blow it down to
a smooth point, with a morphism
$\nu:Y'\to Z$. Let us understand what is $Z$. We have $\rk Pic(Z)=1$.
 Let
$\Lambda\subset|\O_{\P^3}(1)|$ be the pencil of hyperplanes
through $l$, $\f:=\nu\circ\mu\circ\epsilon$, and $\f_*\Lambda=\A$. Then we have
$$-K_Y=\epsilon^*(-K_X)-5E=4H+(12-5)E$$
moreover $A:=(\nu\circ\mu)_*E\in\A$
and this yields
$$-K_Z\sim 7A$$
By construction  $6E$ is a Cartier divisor. Therefore we find a Fano 3-fold of index greater than 1 and it is easy to
realize that $Z\iso\P(1,1,2,3)$, and $A\sim\O(1)$.
Via the
map $\f$ we
can also easily understand the elements of $|\O_Z(6)|$. Let $F\subset|\O_{\P^3}(3)|$ be a cubic surface with $\mult_E(F)=3$. Then
$$-K_Y=\epsilon^*(F+H)-5E=F_Y+H_Y+E$$
hence we get
$$7A=-K_Z=(\nu\circ\mu)_*(F_Y+H_Y+E)=F_V+A,$$
where $F_V=\f_*(F)\in|\O_Z(6)|$.
 Equivalently
 $S_T\in|\O_Z(6)|$ is birational to a
cubic and hence to a plane.
The pair in b) projects  birationally, from the
singular point, to the
pair in a).
\end{case}
\begin{case}[c, d] The pair $(\P(1,1,1,2),\O(2))$
  is the cone over the Veronese surface,
  a minimal degree 3-fold, with an hyperplane section. It projects
  to $(\P^3,\O(2))$ from two general
  points in $S_T$. The pair
  $(\P(1,1,1,2),\O(4))$ has been
  discussed in Example
  \ref{ex:noether}. The pairs in d)
  projects birationally, from the
  singular point, to  $(\P(1,1,1,2),\O(4))$.
\end{case}

\begin{case}[e, f, g] \label{case:3} The  quadrics in
  $\P^3$ are easy, think for instance to
  $T_{2,3}$. Cubics with isolated
  singularities have been treated in
  Construction \ref{con:33}. Cubics with
  non isolated singularities are
  singular along a line and we can apply
  the proof of case \ref{case:iv} in
  Lemma \ref{lem:fiber}. The pair in f)
  projects, birationally, from the
  singular point to
  $(\P^3,\O(3))$. Similarly all the
  pairs in g) project, from a
  suitable number of points in $S_T$, to
  a pair in e).
\end{case}
\end{proof}

\end{document}